\documentclass[12pt]{article}

\usepackage{amsmath}
\usepackage{multirow} 

\usepackage[cp850]{inputenc}
\usepackage{latexsym,amsfonts}

\usepackage{graphicx}
\usepackage{subcaption}
\usepackage{pdflscape}

\usepackage{caption}

\usepackage{tikz-cd} 
\usepackage{tikz}

\usepackage{enumitem} 
\usepackage{authblk}

\newtheorem{Thm}{Theorem}[section]
\newtheorem{Note}{Note}

\newtheorem{Con}[Thm]{Conjecture}

\newtheorem{Def}[Thm]{Definition}
\newtheorem{cons}[Thm]{Construction}

\author[1]{Mustafa Gezek}
\author[2]{Vladimir D. Tonchev\footnote{Corresponding author, tonchev@mtu.edu}}
\affil[1]{Department of Mathematics, Tekirdag Namik Kemal University, 
Tekirdag, Turkey 59030, mgezek@nku.edu.tr}

\title{\bf On partial geometries arising from maximal arcs}
\affil[2]{Department of Mathematical Sciences, Michigan Technological 
University, Houghton, MI USA 49931, tonchev@mtu.edu}

\begin{document}
\maketitle

\begin{abstract}
The subject of this paper are partial geometries $pg(s,t,\alpha)$
with parameters $s=d(d'-1), \ t=d'(d-1), \ \alpha=(d-1)(d'-1)$,
$d, d' \ge 2$. In all known examples, $q=dd'$ is a power of 2
and the partial geometry arises from a maximal arc of degree $d$ or $d'$
in a projective plane of order $q$  via a known construction
due to Thas \cite{Thas73} and Wallis \cite{W},
with a single known exception
of a partial geometry $pg(4,6,3)$ found by Mathon \cite{Math}
that is not associated
with a maximal arc in the projective plane of order 8.
A parallel class of lines is a set of pairwise disjoint lines that covers 
the point set. Two parallel classes are called orthogonal
if they share exactly one line.
An upper bound on the maximum number of pairwise orthogonal
parallel classes in a partial geometry $G$
with parameters $pg(d(d'-1),d'(d-1),(d-1)(d'-1))$
is proved, and it is shown that
a necessary and sufficient condition for $G$ to arise from a maximal arc of
degree $d$ or $d'$ in a projective plane of order $q=dd'$ is that
 both $G$ and its dual geometry contain sets
of pairwise orthogonal parallel classes that meet the upper bound.
An alternative construction of Mathon's partial geometry
is presented, and the new necessary condition 
is used to demonstrate why this partial geometry
is not associated with any maximal arc 
in the projective plane of order 8. 
 The partial geometries associated
with all known maximal arcs in projective planes of order 16
are classified up to isomorphism, and their parallel classes of
lines and the 2-rank of their incidence matrices are computed.
Based on these results,
some open problems and conjectures are formulated.

\medskip


\end{abstract}

\vspace{5mm}

\section{Introduction}

We assume familiarity with basic facts and notions from
combinatorial design theory  \cite{BJL}, \cite{T88}.

A {\it partial geometry} with parameters $s, t, \alpha$,
or shorty,  $pg(s,t,\alpha)$, is a pair $(P,L)$ of a set $P$
of {\it points} and a set $L$ of {\it lines}, with an incidence relation
between points and lines, satisfying the following axioms:
\begin{enumerate}
\item A pair of distinct points is not incident with more than one
line.
\item Every line  is incident with exactly $s+1$ points ($s\ge 1$).
\item Every point is incident with exactly $t+1$ lines ($t\ge 1$).
\item For every point $p$ not incident with a line $l$, 
there are exactly $\alpha$ lines ($\alpha \ge 1$) 
which are incident with $p$, 
and also incident with some point incident with $l$.
\end{enumerate} 

In terms of the parameters
$s$, $t$, $\alpha$, the number $v=|P|$ of  points, and the number
$b=|L|$ of lines of a partial geometry $pg(s,t,\alpha)$
are given by eq. (\ref{eqvb}).
\begin{equation}
\label{eqvb}
v=\frac{(s+1)(st + \alpha)}{\alpha}, \ b=\frac{(t+1)(st + \alpha)}{\alpha}. 
\end{equation}
If $G=(P,L)$ is a partial geometry $pg(s,t,\alpha)$,
the incidence structure $G'$ having as points the lines of $G$, 
and having as lines the points of $G$, where a point and a line
are incident in $G'$ if and only if the corresponding line and a point 
of $G$ are incident, is a partial geometry $pg(t,s,\alpha)$, called the 
{\it dual} of $G$. 

Partial geometries were introduced by R. C. Bose \cite{Bose}.
In the original Bose's notation, the 
number $t+1$ of lines incident with a point 
is denoted by $r$, and the number $s+1$ of points incident 
with a line is denoted by $k$. 
The $(s,t,\alpha)$-notation was adopted later to match
the notation for generalized quadrangles, which are partial
geometries with $\alpha = 1$ \cite{PT}.
A partial geometry $pg(s,t,\alpha)$ is called {\it proper} 
if $1<\alpha< min(s,t)$.

In this paper, we consider proper partial geometries with parameters 
\begin{equation}
\label{eq1}
s = q-d, \ t = q(d-1)/d, \ \alpha = (q-d)(d-1)/d,
\end{equation}
where $q > d\ge 2$ are integers and $d$ divides $q$.
All known partial geometries with parameters (\ref{eq1}) 
arise from maximal arcs\footnote{A maximal arc of degree $d$ in a projective
plane of order $q=dd'$ is a set $\cal{A}$ of $qd-q+d$ points such that
every line is either disjoint from $\cal{A}$ or meets $\cal{A}$ in exactly $d$ points.}
 of degree $d=2^m$ in projective planes 
of order $q=2^h$, $h\ge 2$,
via a construction due to Thas \cite{Thas73}, \cite{Thas74},
\cite[Construction 41.21]{Thas07}, 
and Wallis \cite{W},
with a single exception for $h=3$ and $m=1$, when in addition 
to the partial geometry $pg(6,4,3)$ arising from a maximal arc 
of degree 2, that is, a hyperoval in $PG(2,8)$, there exists a second 
nonisomorphic partial geometry with the same
parameters that was found by Rudi Mathon \cite{M}, 
and is not associated with any hyperoval in $PG(2,8)$.

If $d=2^m$ and $q=2^h$, ($m < h$), the parameters (\ref{eq1}) can 
be written as
\begin{equation}
\label{eq2}
s=2^h - 2^m, t=2^h - 2^{h-m}, \alpha = (2^m -1)(2^{h-m} -1), (1\le m<h). 
\end{equation}
The parameters (\ref{eq2}) correspond to partial geometries
of Type 1 in the classification of the known proper partial
geometries given in \cite[Theorem 41.31]{Thas07}.
(The only known "improper" partial geometry with parameters (\ref{eq2})
is a generalized quadrangle $pg(2,2,1)$,
($m=1$, $h=2$),
arising from a hyperoval in the projective plane of order 4.)

In Section \ref{S2}, we prove a necessary and sufficient condition
for a partial geometry with parameters (\ref{eq1})
to be associated with a maximal arc of degree $d$ in a projective 
plane of order $q$ in terms of parallel classes of lines
 (Theorem \ref{t2}). 

In Section \ref{S3},
we give an alternative construction of Mathon's partial geometry 
with parameters $s=6$, $t=4$, $\alpha=3$, and 
use Theorem \ref{t2} to show why this partial geometry
is not associated with a maximal arc in $PG(2,8)$.

In Section \ref{S4}, we examine the  partial geometries
associated with maximal arcs in projective planes of
even order $q\le 16$. The partial geometries associated
with all known maximal arcs in projective planes of order 16
are classified up to isomorphism, and their parallel classes of 
lines and the 2-rank of their incidence matrices are computed.
Based on the results of these computations,
some open problems and two conjectures are formulated.

\section{Maximal arcs and partial geometries}
\label{S2}

Let $k$ and $d$ be positive integers such that $k > d > 1$.
A $(k,d)$-arc (or an arc of size $k$ and degree $d$) 
in a projective plane $\cal{P}$ of order $q$ is 
a set $\cal{A}$ of $k$ 
points such that $d$ is the greatest number of 
collinear points in $\cal{A}$. It follows that 
$k \le dq-q+d$, 
and the equality $k=dq-q+d$ holds if and only
if every line is either disjoint from $\cal{A}$
 or intersects $\cal{A}$ in exactly $d$ points. 
A $(dq-q+d,d)$-arc is called a {\it maximal arc}.
A maximal arc of degree $d=2$ is  called also a {\it hyperoval}.

A necessary condition for the existence
of a maximal arc of degree $d$ in a projective
plane of order $q$ is that $d$ divides $q$.
If $\cal{A}$ is a maximal $(dq-q+d,d)$-arc in a projective
plane $\cal{P}$ of order $q$, then the lines of $\cal{P}$ 
that are disjoint from $\cal{A}$ form a maximal 
$(q(q-d+1)/d,q/d)$-arc $\cal{A}'$ in the dual plane 
$\cal{P}'$, called the {\it dual arc} of $\cal{A}$. 

Maximal arcs of degree $d$,  $1< d < q$, do
not exist in any Desarguesian plane of odd order
$q$ (Ball, Blokhuis, and Mazzocca \cite{BBM}),
as well as in any of the four projective planes of order 9
(Lunelli and Sce \cite{LS}),
and are known to exist in every Desarguesian plane of
order $q=2^t$ (Denniston \cite{Den}, Hamilton 
and Mathon \cite{HM1}, \cite{HM2},
Mathon \cite{Math}, Thas \cite{Thas74}),
 as well as in some non-Desarguesian planes of even order
(Gezek, Mathon and Tonchev \cite{GMT},  Gezek, Tonchev and Wagner 
\cite{GTW}, Hamilton  \cite{H1}, \cite{H2}, \cite{H3},
Hamilton, Stoichev and Tonchev \cite{HST},
Penttila, Royle, and Simpson \cite{PRS},
Thas \cite{Thas80}).

If $\cal{P}$ is a projective plane of order
 $q=dd'$, then a maximal arc $\cal{A}$ in  $\cal{P}$
 of degree $d$ has 
$d(q-d'+1)$ points, while its dual arc $\cal{A}'$
is of degree $d'$ and has $d'(q-d+1)$ points.
The set of points of a maximal $(d(q-d'+1),d)$-arc $\cal{A}$
and the non-empty intersections of $\cal{A}$ with lines of $\cal{P}$
considered as blocks, form a Steiner 2-$(d(q-d'+1),d,1)$ design $D$.
We say that $D$ is associated with the maximal arc $\cal{A}$,
or that $D$ is embeddable in $\cal{P}$ as a maximal arc.
Similarly, the points of the dual arc $\cal{A}'$ define
a Steiner 2-$(d'(q-d+1),d',1)$ design $D'$
embeddable in the dual plane $\cal{P}'$.

A necessary and sufficient condition for a Steiner 
2-$(d(q-d'+1),d,1)$ design $D$ to be embeddable as a maximal
$(d(q-d'+1),d)$-arc in a projective plane of order $q=dd'$ was proved
by the second author in \cite{ton-res}.
This condition was used in \cite{GTW} to show that
five of the known projective planes of order 16 contain 
maximal $(52,4)$-arcs
whose associated Steiner 2-$(52,4,1)$ designs are embeddable in 
two nonisomorphic planes of order 16.
\begin{cons}
\label{con}
(Thas \cite{Thas73}, \cite{Thas07}, Wallis \cite{W}). 
Let $\cal{A}$ be a maximal $(qd-q+d,d)$-arc in a projective plane
$\cal{P}$ of order $q=dd'$, $2 \le d < q$, and let $X$ be the
set of points of $\cal{P}$.
Let $P = X \setminus \cal{A}$, and let $L$ be the set of all lines 
of $\cal{P}$ that intersect $\cal{A}$ in $d$ points. Then $G=(P,L)$ 
is a partial geometry $pg(s,t,\alpha)$ with parameters (\ref{eq1}).
Similarly, the dual arc $\cal{A}'$ determines a partial geometry $G'$ 
with parameters $t,s,\alpha$, being the dual geometry of $G$.
\end{cons}
We note that if $q=dd'$, then
the parameters (\ref{eq1}) of a partial geometry associated
with a maximal arc $\cal{A}$ of degree $d$
via Construction \ref{con}  can be rewritten as in eq. (\ref{eq3}),
while the numbers of points and lines (\ref{eqvb}) can be written
as in eq. (\ref{eqnvb}).
\begin{equation}
\label{eq3}
s = d(d' -1), \ t = d'(d-1), \ \alpha = (d-1)(d' -1). 
\end{equation}
\begin{equation}
\label{eqnvb}
v=(s +1)(dd' +1), \ b=(t +1)(dd' +1). 
\end{equation}
Let $d\ge 2$, $d' \ge 2$ be integer numbers,
 and let $G=(P,L)$ be a partial geometry   
with parameters $s$, $t$, $\alpha$ given by eq. (\ref{eq3}).
 Thus,  by eq. (\ref{eqnvb}), the line size $s+1$ of $G$  divides 
the number of points $v$, while $t+1$,
that is, the line size the dual geometry $G'$, divides $b$. 
A {\it parallel class} of $G$
 is a set of $v/(s+1)=dd'+1$ pairwise disjoint lines  of $G$.
Similarly, a parallel class of the dual geometry $G'$
is a set of $b/(t+1)=dd' +1$
pairwise disjoint lines of $G'$.

\begin{Def}
Two parallel classes of $G$ (resp. $G'$) 
are called {\it orthogonal} if 
they share exactly one line.
\end{Def}

The following theorem gives an upper bound on the maximum number
of pairwise orthogonal parallel classes.

\begin{Thm}
\label{t1}
 Let $G=(P,L)$ be a partial geometry $pg(s,t,\alpha)$
 with parameters (\ref{eq3}), where $d\ge 2$, $d'\ge 2$ 
are given integer numbers, and let $G'$ be its dual partial geometry.\\
 (a) If $C$ is a set of $m$ pairwise orthogonal parallel classes
 of $G$ then  
\begin{equation}
\label{eq6}
m \le d(dd' - d' +1),
\end{equation} 
and the equality $m = d(dd' - d' +1)$ holds if and only if every 
line of $G$ appears in   exactly $d$ parallel classes from $C$.
 
(b) If $C'$ is a set of $m'$ pairwise orthogonal parallel classes
 of $G'$ then
 \begin{equation}
\label{eq6'}
m' \le d'(d'd -d +1),
\end{equation}
and the equality $m' = d'(d'd - d +1)$ holds if and only if every 
line of $G'$ appears in   exactly $d'$ parallel classes from $C'$.
\end{Thm}
{\bf Proof}. (a) Let $k_i$ denote the number of parallel classes from
$C$ that contain  the $i$th line $L_i$ of $G$, $1\le i \le b$,
where the number $b$ of lines of $G$ is
equal to $b=(d'(d -1)  +1)(dd' +1)$ by eq.
 (\ref{eqnvb}) and (\ref{eq3}). 
Let $C_1,  \ldots, C_m$ be the parallel classes from $C$.
Counting in two ways the ordered pairs $(C_i, L_t)$, where
$L_t$ belongs to $C_i$ gives
\begin{equation}
\label{eq7}
\sum_{i=1}^b k_i = m(dd' +1).
\end{equation}
Similarly, counting in two ways the ordered pairs
$(\{ C_i, C_j \}, L_t)$, where $L_t$ belongs to $C_i$ and $C_j$
($i \neq j$), gives
\begin{equation}
\label{eq8}
\sum_{i=1}^b k_{i}(k_i - 1) = m(m-1).
\end{equation}
Adding equations (\ref{eq7}) and (\ref{eq8}) gives
\[ \sum_{i=1}^b k_{i}^2 = m(m + dd').
\]
Applying the Cauchy-Schwarz inequality, we have
\[ (\sum_{i=1}^b  k_i )^2 =m^{2}(dd' +1)^2 \le b\sum_{i=1}^b k_{i}^2 =
(d'(d-1)+1)(dd'+1)m(m+dd'), \]
whence
\[ m(dd' +1) \le (d'(d-1)+1)(m+dd'). \]
Solving the last inequality for $m$ implies (\ref{eq6}).
Clearly,  the equality 
\[ m=d(dd' -d' +1) \]
holds if and only if
\[ k_1 = k_2 = \cdots = k_b =\frac{1}{b}\sum_{i=1}^b k_i =
\frac{d(dd' -d' +1)(dd' +1)}{(d'd -d' +1)(dd' +1)} =d. \]
This completes the proof of part (a). 

(b) Switch $d$ and $d'$, $b$ and $v$, $s$ and $t$
in the proof of part (a).     $\Box$

The following theorem shows that a partial geometry $G$ with
parameters (\ref{eq3}) arises from a maximal arc
if and only if $G$ and its dual geometry $G'$ both meet the bounds
of Theorem \ref{t1}.
\begin{Thm}
\label{t2}
Let $d\ge 2$, $d' \ge 2$ be integer numbers, 
and let $G$ 
be a partial geometry  $pg(s,t,\alpha)$  
with parameters  (\ref{eq3}).

A necessary and sufficient condition for $G$
to be associated with a maximal arc of degree $d$ in  
a projective plane $\cal{P}$ of order
$q=dd'$ is that the following two conditions hold:

(i) $G$ 
 admits a set of  pairwise orthogonal parallel classes
 that meets the upper bound (\ref{eq6}) of Theorem \ref{t1}, part (a).

(ii) The dual geometry $G'$
admits a set of  pairwise orthogonal parallel classes
 that meets the upper bound (\ref{eq6'}) of Theorem \ref{t1},
part (b).
 \end{Thm}

{\bf Proof.} 
Assume that $\cal{P}$ is a projective plane of order  $q=dd'$ 
with a maximal $(qd-q+d,d)$-arc $\cal{A}$, 
and let $\cal{A}'$ be the dual $(qd'-q+d',d')$-arc.
Let $G$ be the partial geometry with parameters given by (\ref{eq3})
and (\ref{eqnvb}) that arises from $\cal{A}$ via Construction \ref{con}.
Thus, if $X$ is the point set of $\cal{P}$ then $X\setminus \cal{A}$
is the set of points of $G$.

A {\it pencil} is a set of $q+1$ lines of $\cal{P}$ that pass through a
point $x$.
If $x$ is a point from the arc $\cal{A}$, the pencil through $x$ determines 
\[ v/(s+1)=q+1 =dd' +1 \]
pairwise disjoint lines of $G$ via Construction \ref{con}, 
and these $q+1$ lines of $G$ form a parallel class.
Since every two points $x, y \in \cal{A}$, $x\neq y$ are incident 
with exactly one line $l$ of $\cal{P}$, the pencils through $x$ and $y$
determine two parallel classes $C_1$, $C_2$ of lines of $G$ that share 
exactly one line, being the restriction of $l$ 
on the point set $X\setminus \cal{A}$ of $G$. 
Thus $C_1$ and $C_2$ are orthogonal.
It follows that the $qd-q+d=d(dd' -d' +1)$ points of $\cal{A}$
determine a set of $d(dd' -d' +1)$ pairwise orthogonal classes of lines of $G$
that meets the bound (\ref{eq6}) of Theorem \ref{t1}.

Similarly, the $d'(d'd-d+1)$ points of $\cal{A}'$ determine
a set of $d'(d'd-d+1)$ pairwise orthogonal parallel  classes of lines 
of the dual geometry $G'$
that meets the bound (\ref{eq6'}) of Theorem \ref{t1}.  

Suppose now that $G =(P,L)$ is a partial geometry with parameters
(\ref{eq3}) for some integers $d\ge 2$, $d' \ge 2$,
and assume that $G$ and its dual geometry $G'$ satisfy the 
conditions (i) and (ii) respectively. Then one can construct 
a projective plane $\cal{P}$ of order $q=dd'$ as follows.

Let $C=\{ C_1,\ldots,C_m \}$ be a set of $m$ pairwise 
orthogonal parallel 
classes of lines of $G$, where $m=d(dd'-d'+1)$
meets the bound (\ref{eq6}),  
and let $C' = \{ C'_1,\ldots, C'_{m'} \}$ be  a set of 
$m'$ pairwise orthogonal parallel classes of lines of
the dual geometry $G'$, where $m' = d'(d'd-d+1)$
meets the bound (\ref{eq6'}). 
According to Theorem \ref{t1},
every line of $G$ appears in exactly $d$ parallel
classes from $C$, and every line of $G'$ appears in 
exactly $d'$ parallel classes from $C'$.
By (\ref{eq3}) and (\ref{eqnvb}), every parallel class 
$C_i \in C$ , as well as every parallel class $C'_i \in C'$ ,
consists of $q+1$ pairwise disjoint lines.

We define an incidence structure $\cal{P}$
with a set of points $X$ and a collection of lines $\cal{L}$,
where $X$ consists of the $v=|P|$ points of $G$ plus
 $m=d(dd'-d'+1)$
new points labeled by the $m$ parallel classes from $C$,
and $\cal{L}$ consists of $b=|L|$ lines labeled by the lines
of $G$ plus $m' = d'(d'd-d+1)$ lines labeled by the $m'$
parallel classes from $C'$.
Thus, by (\ref{eq3}) and (\ref{eqnvb}), $\cal{P}$ has
$q^2 + q + 1$ points and $q^2 + q + 1$ lines.

A line $l^*$ of $\cal{P}$ which is labeled by a line $l$ of $G$
consists of the $s+1$ points of $l$ and $d$ points labeled
by the $d$ parallel classes from $C$ that contain
the line $l$. It follows from (\ref{eq3}) that $l^*$ 
is incident with $s+1 + d =  q+1$ points.

A line $l'_i$ of $\cal{P}$ which is labeled
by a parallel class $C'_i \in C'$,
consists of the $q+1$ points of $G$ that correspond
to the $q+1$  lines of $G'$ belonging to $C'_i$,
$1\le i \le m'$.
It follows from (\ref{eq3}) that every point of $\cal{P}$
which is also a point of $G$ is incident with
$t+1 + d' = q+1$ lines of $\cal{P}$,
and every point of $\cal{P}$ which is labeled
by a parallel class $C_i$ of $G$ is also incident
with $q+1$ lines of $\cal{P}$.

Thus, $\cal{P}$ = $(X, \cal{L})$ is an incidence structure
with $|X|$=$q^2 +q + 1$ points and $|\cal{L}|$=$q^2 + q + 1$ lines,
such that every line is incident with $q+1$ points, and every
point is incident with $q+1$ lines. To show that $\cal{P}$
is a projective plane of order $q$, it is sufficient to
check that every two lines of $\cal{L}$ meet in exactly one
point of $X$.

Any two distinct parallel classes  $C'_i, C'_j \in C'$ 
share exactly one line $l'$ of $G'$
due to the orthogonality condition.
Consequently, the two lines of $\cal{P}$ labeled by
$C'_i$ and  $C'_j$ 
 meet in
exactly one point $x\in P$, being the point of $G$ that
corresponds to the line $l'$ of $G'$.

A line $l^*$ of $\cal{P}$ which is labeled
by a line $l$ of $G$ meets any line of $\cal{P}$
which is labeled by a parallel class
$C'_i \in C'$ in exactly one point, because
$l$ is incident with exactly one of the $q+1$
points  of $G$ that correspond to the $q+1$ parallel
lines of $G'$ belonging to $C'_i$.

Let $l^*_1$, $l^*_2$ be lines of $\cal{P}$ that are
labeled by two distinct lines $l_1$, $l_2$ of $G$.

If  $l_1$ and $l_2$ belong to a parallel class
$C_i \in C$, then
$l^*_1$ and $l^*_2$ meet in exactly one point,
being the point of $\cal{P}$  labeled by $C_i$.
On the other hand, if $l_1$ and $l_2$ are two distinct lines of $G$
that belong to two different parallel classes from $C$, then
$l^*_1$ and $l^*_2$ cannot share any point labeled by a parallel
class from $C$, and can possibly share at most one point, being 
a point of $G$.

To prove that every two lines of $\cal{P}$ share a point, 
we count the set of ordered pairs
\[
S =\{ (x, \{ l^*_i, l^*_j \}) \},
\] 
where $x\in X$ is a point of $\cal{P}$, and $l^*_i, l^*_j \in\cal{L}$,
$1\le i < j \le q^2 + q + 1$,
 are two
distinct lines of $\cal{P}$ that are both incident with $x$.

Since for every two distinct lines $l^*_i$, $l^*_j$
there is at most one point $x$ incident with both $l^*_i$
and $l^*_j$,  we have
\begin{equation}
\label{xl1l2}
|S| \le { q^2 + q + 1 \choose 2}.
\end{equation} 

On the other hand, since every point $x$ belongs to $q+1$ lines, 
there are
${ q+1 \choose 2}$ pairs of lines that are incident with $x$.
Thus, counting the ordered pairs from $S$ by ranging $x$ over the
set of all $q^2 + q + 1$ points  implies that
$|S|$ is given by eq. (\ref{x}). 
\begin{equation}
\label{x}
 |S|=(q^2 + q + 1){q+1 \choose 2}=\frac{(q^2+q +1)(q+1)q}{2}.
\end{equation}
Since
\begin{equation}
\label{bin}
 \frac{(q^2+q +1)(q+1)q}{2} = {q^2 + q + 1 \choose 2},
\end{equation} 
it follows from (\ref{xl1l2}), (\ref{x}) and (\ref{bin})
that each pair of distinct lines share exactly one point,
thus $\cal{P}$ is a projective plane of order $q=dd'$.

Since every line of $\cal{P}$ which is labeled by a line of $G$
meets the set $\cal{A}$ of $m=d(dd'-d'+1)$ points labeled by the
parallel classes from $C$ in exactly $d$ points, and every line
of $\cal{P}$ which is labeled by a parallel class from $C'$
is disjoint from $\cal{A}$, the set $\cal{A}$ is a maximal
$(d(dd'-d'+1),d)$-arc in $\cal{P}$. Similarly, the
$m'=d'(d'd-d+1)$ lines of $\cal{P}$ labeled by the parallel 
classes of $C'$,
determine a maximal $(d'(d'd-d+1),d')$-arc $\cal{A}'$ in the 
dual plane
of $\cal{P}$, and $\cal{A}'$ is the dual arc of $\cal{A}$. $\Box$. 

\begin{Note}
\label{n1}
{\rm
In all  known partial geometries with parameters of the form 
(\ref{eq3}), $d$ and $d'$ are both powers  of 2. 
It is an interesting 
open question whether partial geometries for other values 
of $d$ and $d'$ may exist. 

For example, it is not known if a partial geometry $pg(6,6,2)$,
($d=d'=3$) exists or not.
However, a partial geometry with these parameters
cannot be associated with a maximal $(21,3)$-arc in a projective plane 
of order 9, because no such maximal arcs exist in any of the
four projective planes of order 9 \cite{LS}.

Another open small parameter set is $s=12, \ t=10, \ \alpha =8$, 
that corresponds to $d=3$ and $d'=5$. The existence
of a partial geometry $pg(12,10,8)$ 
that satisfies the conditions 
of Theorem \ref{t2} would imply the existence of a projective
plane of order 15 with maximal arcs of degree 3 and 5.
}
\end{Note} 

\section{Partial geometries $pg(4,6,3)$}
\label{S3}

A {\it strongly regular graph} $\Gamma$ with parameters $n$, $k$, $\lambda$,
$\mu$ (or srg$(n,k,\lambda,\mu)$) is an undirected graph 
without loops or multiple edges, having the following properties:
\begin{itemize}
\item $\Gamma$ has $n$ vertices.
\item Every vertex has exactly $k$ neighbors.
\item Every two adjacent vertices have exactly $\lambda$ common neighbors.
\item Every two nonadjacent vertices have exactly $\mu$ common neighbors.
\end{itemize}

The complementary graph $\bar{\Gamma}$ of a strongly regular graph
$\Gamma$ with parameters $n, k, \lambda, \mu$
is also strongly regular, 
with parameters $\bar{n}, \bar{k}, \bar{\lambda}, \bar{\mu}$
given by (\ref{eqC}).
\begin{equation}
\label{eqC}
 \bar{n}=n, \ \bar{k} = n-1-k, \ \bar{\lambda}=n - 2k + \mu -2, \ \bar{\mu}
= n - 2k +\lambda.
\end{equation}

Strongly regular graphs were introduced by Bose in the same paper
\cite{Bose},  where he introduced partial geometries.

The eigenvalues of the (0,1)-adjacency matrix of a strongly
regular graph $\Gamma$ with parameters $n, k, \lambda, \mu$,
as well as their multiplicities, are easily expressed in terms of the
graph parameters \cite{Bose}: $k$ is a simple eigenvalue,
and up to multiplicity, there are two more
eigenvalues $\rho_1, \rho_2$,
being the solutions of the
quadratic equation (\ref{eqr1r2}).
\begin{equation}
\label{eqr1r2}
x^2 +(\mu -\lambda)x +\mu -k=0.
\end{equation}
The following statement is a special case of a more general
result due to Hoffman \cite{Hoff} that applies to
regular graphs.
\begin{Thm}
\label{hof}
 (Hoffman bound)
Let $\Gamma$ be a strongly regular graph with parameters $n,k,\lambda,\mu$,
and 
let $\rho$ be the smallest eigenvalue of the (0,1)-adjacency matrix of
$\Gamma$,
being the negative root of equation (\ref{eqr1r2}).
The size of any coclique $C$ of $\Gamma$ satisfies the inequality
\begin{equation}
\label{hofb}
c=|C|\le \frac{n(-\rho)}{k-\rho},
\end{equation}
and the equality holds if and only if every vertex outside $C$ is adjacent to
exactly
\[ d=\frac{kc}{n-c} \]
vertices of $C$.
\end{Thm}

If $G=(P,L)$ is a partial geometry $pg(s,t,\alpha)$
 with point set $P$ and line set $L$,
the {\it point graph} $\Gamma_P$ of $G$ is the graph with vertex set $P$,
where two vertices are adjacent if the corresponding points of $G$ are
collinear. The {\it line graph} $\Gamma_L$ of $G$
is the graph having as vertices the lines of $G$, where two lines are adjacent
if they share a point. Both $\Gamma_P$ and $\Gamma_L$ are strongly regular 
graphs \cite{Bose}. The parameters $n, k, \lambda, \mu$ of $\Gamma_P$
are expressed in terms of $s, t, \alpha$ as in eq. (\ref{eqP}),
while the parameters $n', k', \lambda', \mu'$ of $\Gamma_L$
are given by eq. (\ref{eqL}).
\begin{equation}
\label{eqP}
n=(s+1)(st+\alpha)/\alpha,  k = s(t+1),  \lambda = s-1 +t(\alpha -1), \mu =
\alpha(t+1). 
\end{equation}
\begin{equation}
\label{eqL}
n' = (t+1)(st+\alpha)/\alpha,  k' =t(s+1),  \lambda' = t-1 +s(\alpha -1),
 \mu' = \alpha(s+1). 
\end{equation} 

A strongly regular graph $\Gamma$ whose parameters
$n, k, \lambda, \mu$ can be written as in eq. (\ref{eqP}) for some
integers $s, t, \alpha$ is called {\it pseudo-geometric},
and  $\Gamma$ is called {\it geometric}
 if there exists a partial geometry $G$
with parameters $s, t, \alpha$ such that $\Gamma$ is the point graph 
of $G$; otherwise  $\Gamma$ is  {\it non-geometric}. 

Applying the Hoffman bound (\ref{hofb}) from Theorem \ref{hof} 
to the complementary graph
$\bar{\Gamma}$ of a pseudo-geometric graph $\Gamma$ with parameters
(\ref{eqP}) implies that every clique of $\Gamma$ is
of size smaller than or    
equal to $s+1$, and every clique $C$ of maximum size $s+1$ has the
property that every vertex outside $C$ is adjacent to exactly
$\alpha$ vertices from $C$. 
\begin{Thm}
\label{bt}
(Bose \cite{Bose}).
A pseudo-geometric strongly regular graph $\Gamma$
with parameters (\ref{eqP}) is geometric if and only if
$\Gamma$ possesses a set of $b=(t+1)(st+\alpha)/\alpha$ cliques
of size $s+1$,
 every two of which share at most one vertex.
\end{Thm}
Next, we consider  partial geometries 
with parameters  $s, t, \alpha$ of the form (\ref{eq3}) for some integer
$d\ge 2$ and $d' = 2$, namely 
\begin{equation}
\label{nsd}
s=d, \ t=2d-2, \ \alpha = d-1.
\end{equation}
The parameters (\ref{eqP})
of the point graph of a partial geometry with parameters
(\ref{nsd}) are given by eq. (\ref{srgd2}).
\begin{equation}
\label{srgd2}
n = (d+1)(2d+1),  k = d(2d-1),  \lambda = (d-1)(2d-3),  \mu = (d-1)(2d-1). 
\end{equation}
The {\it triangular graph} $T(m)$, where $m\ge 4$ is an integer, 
has as vertices the unordered 2-subsets of $\{ 1, 2,\ldots, m \}$, 
where two distinct 2-subsets
are adjacent in $T(m)$ whenever they are not disjoint. 
The graph $T(m)$ is strongly regular with parameters
 (\ref{Tm}),
while its complementary graph $\bar{T}(m)$ has parameters (\ref{eqTbar}).
\begin{equation}
\label{Tm}
 n={ m \choose 2}, \ k = 2(m-2), \ \lambda = m-2, \ d = 4.
\end{equation}
\begin{equation}
\label{eqTbar}
\bar{n} = { m \choose 2},  \bar{k} = { m-2 \choose 2}, \bar{\lambda} =
 { m-4 \choose 2},  \bar{\mu} = { m - 3 \choose 2}.
\end{equation}
Let $d \ge 2$ be an integer. By
(\ref{eqTbar}), the parameters of $\bar{T}(2d+2)$ are
\begin{equation}
\label{m+2}
\bar{n} = (d+1)(2d+1),  \bar{k} = d(2d-1),  \bar{\lambda} = (d-1)(2d-3),  
 \bar{\mu} = (d-1)(2d-1). 
\end{equation}
Since the parameters (\ref{m+2}) and (\ref{srgd2})
coincide,  the graph $\bar{T}(2d+2)$
is pseudo-geometric for $s, t, \alpha$ given by (\ref{nsd}).
It is known that every strongly regular graph with parameters 
(\ref{Tm}) is isomorphic to $T(m)$, except when $m=8$, in which case
in addition to $T(8)$, there are three other graphs \cite{Chang},
\cite{Hoffman},
known as the Chang graphs.

The question about the values of $d$ for which a strongly regular
graph with parameters (\ref{m+2}) is geometric has been settled
in the following cases:
\begin{enumerate}
\item $d=2^i$, $i\ge 0$. The graph $\bar{T}(2^{i+1}+2)$ is geometric,
and the corresponding partial geometry  
arises from a hyperoval in a projective plane of order $2^{i+1}$
 via Construction \ref{con}.
\item $d=3$.  Neither $\bar{T}(8)$, nor any of the Chang graphs
 is geometric, thus, a partial geometry with parameters
$s=3, \ t=4, \alpha =2$ does not exist (De Clerck \cite{DeCl}).
\item $d=4$. The graph $\bar{T}(10)$ is geometric, and up to isomorphism,
there exist exactly two partial geometries $pg(4,6,3)$ (Mathon \cite{M}).
The dual geometry of
one of these two geometries arises from a hyperoval
in the projective plane of order 8, $PG(2,8)$, via Construction \ref{con},
while the dual of the second geometry is not associated with any hyperoval
in $PG(2,8)$.
\item $d=5$. The graph $\bar{T}(12)$ is not geometric, 
thus a partial geometry
$pg(5,8,4)$ does not exist (Lam, Thiel, Swiercz, and McKay \cite{Lam}).
\end{enumerate}
We will give a brief description of the two partial geometries
$pg(4,6,3)$ having $\bar{T}(10)$ as point graph. 
A convenient way to describe a partial geometry
$G_1$ associated with a
maximal arc of degree 4 (dual arc of a hyperoval) in $PG(2,8)$ 
is by using a collineation 
$f$ of order 9 acting fixed-point-free on the set of points
$P=\{ 1,2,\ldots, 45 \}$, as
\[ f_P = (1,2,\ldots, 9)(10,\ldots,18)\cdots(37,\ldots,45), \]
and on the set of 63 lines $L = \{ L_1,\ldots, L_{63} \}$, as
\[ f_L =  (L_1,\ldots, L_9)(L_{10},\ldots,L_{18})\cdots (L_{55},\ldots,L_{63}). \]
Representatives of the seven orbits of lines are listed in Table \ref{tabL'}.
The order of the full automorphism group of $G_1$ is 1512,
which is also the order of the stabilizer of a hyperoval in 
$PG(2,8)$.
The 2-rank (that is, the rank over the finite field of order 2)
of the incidence matrix of $G_1$ is equal to 28. We note
that 28 is also the 2-rank of the incidence matrix of the 
projective plane of order 8.

We define a graph $\cal{L}$ having as vertices the lines of $G_1$,
where two lines are adjacent in $\cal{L}$ if they  are disjoint,
or in other words, $\cal{L}$ is the complementary graph of the line 
graph of $G_1$. Clearly, the maximum clique size in $\cal{L}$
is $45/5 =9$, and every clique of size 9 is a parallel 
class of lines. 
Using the clique finding algorithm Cliquer 
developed by Niskanen and  \"Osterg\r ard \cite{cliquer},
one quickly finds by computer that  $\cal{L}$ 
contains exactly 28 cliques of size 9, every two sharing
one vertex, that give a set $C$
of 28 pairwise orthogonal parallel classes of lines,  
listed in Table \ref{pclG1},
where the indices of the lines in each
parallel class are given. Note that 28 is equal to the upper 
bound (\ref{eq6}) of Theorem \ref{t1}, (a).
Similarly, the parallel classes of lines of the dual geometry
$G'_1$, each consisting of nine pairwise disjoint lines
of $G'_1$, can be found as  cliques of size 9
 in the complementary graph of the line graph of $G'_1$,
or equivalently, the complementary graph of the point graph of $G_1$,
the latter being isomorphic to the triangular graph $T(10)$.
Clearly, there are exactly ten such cliques meeting pairwise
in one vertex.
The set $C'$ of ten specific cliques found by Cliquer are listed in
Table \ref{pclPG1}, where each line of the dual geometry
 $G'_1$ is labeled by a point of $G_1$.
We note that  the set of ten pairwise orthogonal parallel 
classes meets the upper bound (\ref{eq6'}) of Theorem \ref{t1} (b).
Thus, a projective plane of order 8 is uniquely
determined from $G_1$ and the sets of parallel classes $C$ and $C'$
by the construction described in the proof of Theorem \ref{t2}.

The second partial geometry $G_2$ with parameters $s=4, \ t=6, \ \alpha=3$,
found by Mathon \cite{M}, can be described in terms of
a  permutation $g$  of order 6
and representatives of the line orbits under the group 
generated by $g$ as follows.
The permutation $g$ acts on the 45 points as
\[ g_P =(1,\ldots,6)(7,\ldots,12) \cdots (31,\ldots,36)(37,38,39)(40,41)(42,43)(44)(45), \]
and $g$ acts on the set of 63 lines $\{ l_1,\ldots,l_{63} \}$
as
\[ g_L = (l_1,\ldots,l_{6}) \cdots (l_{49},\ldots,l_{54})(l_{55},l_{56},l_{57})(l_{58},l_{59},l_{60})(l_{61},l_{62})(l_{63}). \]
We take as line orbit representative the first line from each orbit.
These representatives are given in Table \ref{lMr}.
The order of the full automorphism group of $G_2$ is 216 \cite{M}.
The 2-rank of the incidence matrix of $G_2$ is equal to 34,
thus, $G_2$ does not arise from 
a maximal arc in a projective plane of order 8, 
because the 2-rank of the incidence matrix of 
the (unique up to isomorphism)
projective plane of order 8 is equal to 28.
To find all parallel classes of lines of $G_2$,
we consider each parallel class as a  
clique of size 9 in the complementary graph 
of the line graph of $G_2$.
Using Cliquer, we found that $G_2$ contains only one 
parallel class of lines.
Thus, $G_2$ does not satisfy the first codition of
Theorem \ref{t2}, and consequently, 
this partial geometry is not associated with a maximal arc of degree 4
(dual arc of a hyperoval) in a projective plane of order 8.
 The indices of the lines
form the unique parallel class are
\[ 25,  \ 26,  \ 27,  \ 28,  \  29,  \  30,  \  61,  \  62,  \  63, \]
and the points of these lines are listed
in Table \ref{parlines}.
 The parallel classes of lines of the dual geometry
$G'_2$, each consisting of nine pairwise disjoint lines
of $G'_2$, can be found as maximal cliques of size 9
 in the complementary graph of the point graph of $G_2$,
 being isomorphic to the triangular graph $T(10)$.
There are exactly ten such cliques,
and the cliques found by Cliquer are listed in
Table \ref{pclPG2}, where each line of the dual geometry
 $G'_2$ is labeled by a point of $G_2$.
The  partial geometries $pg(4,6,3)$ having $\bar{T}(10)$ 
as a point graph were classified by Mathon \cite{M}
by enumerating and classifying up 
to isomorphism the collections of 63 5-cliques of 
$\bar{T}(10)$ meeting pairwise in at most one vertex.
A {\it pencil} through a point $x$ of a partial geometry $pg(4,6,3)$
is the set of lines though $x$. To reduce the search,
Mathon \cite{M} enumerated and classified the pencils through 
a pair of vertices of $\bar{T}(10)$, and analyzed their completion  
to a partial geometry.
An alternative construction of Mathon's partial geometry
based on over-large Steiner systems was proposed recently by
Reichard and Woldar \cite{RW}.

Using Cliquer, we enumerated and classified up to isomorphism
 all partial geometries $pg(4,6,3)$ by 
computing all maximal cliques in a graph $\Omega$ having as vertices the
$10!/(2^{5}5!)=945$ 5-cliques of $\bar{T}(10)$,
where two 5-cliques of $\bar{T}(10)$  are adjacent in $\Omega$ 
if and only if they share at most one vertex.
Clearly, any set of 63 5-cliques of $\bar{T}(10)$ that is a clique in $\Omega$
is the line set of a partial geometry $pg(4,6,3)$, and vice versa.
It took Cliquer less than one minute 
on a personal computer MacBook Pro
to find all cliques of size 63 in $\Omega$, an their number is
19200. 

According to the order of  the stabilizer
of a 63-clique in the symmetric group $S_{10}$, which acts as 
the full automorphism group of both $\bar{T}(10)$ and $\Omega$, 
the 19200 63-cliques of $\Omega$ are partitioned into two sets: 
2400 cliques with a stabilizer of order 1512, and the remaining 
16800 63-cliques with a stabilizer of order 216.
Since
\[ \frac{10!}{1512} + \frac{10!}{216} = 2400 + 16800 = 19200, \]
it follows that the 19200 distinct partial geometries $pg(4,6,3)$
having $\bar{T}(10)$ as point graph, are split into two ismoprphism classes,
one geometry with automorphism group of order 1512, being isomorphic to 
the geometry arising from a maximal arc
of degree 4 in $PG(2,8)$, and a second geometry with 
automorphism group of order 216, which is isomorphic to Mathon's geometry. 
Thus, our computations
give an alternative and independent confirmation of Mathon's 
classification \cite{M},
and Theorem \ref{t2} provides an alternative explanation why 
the dual of Mathon's geometry 
does not arise from a hyperoval in a projective plane of order 8.

\begin{table}
\centering
\begin{tabular}{|l|l|}
\hline
$L_1$   & 9 11 22 34 42  \\
\hline
$L_{10}$ &  1  2 18 22 30\\
\hline
$L_{19}$ &  1  3 15 27 42\\
\hline
$L_{28}$ & 11 12 14 19 30\\
\hline
$L_{37}$ & 10 19 40 42 45\\
\hline
$L_{46}$ & 1  5 19 32 41\\
\hline
$L_{55}$ & 19 28 29 34 44\\
\hline
\end{tabular}
\caption{Line orbit representatives of $G_1$}
\label{tabL'}
\end{table}

\begin{table}
\centering
\begin{tabular}{|l|}
\hline
   7 10 13 17 36 41 42 48 62 \\
   2 12 14 17 31 37 45 52 57\\
   3 13 15 18 32 37 38 53 58\\
 5 26 28 33 45 47 48 50 58\\
   2 23 30 34 42 47 53 54 55\\
   6 12 16 18 35 40 41 47 61\\
   3 17 19 22 23 28 44 59 61\\
   2 16 21 22 27 36 43 58 60\\
   6 11 22 25 26 31 38 55 62\\
   9 14 19 20 25 34 41 56 58\\
   9 10 12 15 29 43 44 50 55\\
   9 21 28 32 40 51 52 54 62\\
   4 25 32 36 44 46 47 49 57\\
   6 27 29 34 37 48 49 51 59\\
   7 19 30 35 38 49 50 52 60\\
   8 11 14 18 28 42 43 49 63\\
   7 12 23 26 27 32 39 56 63\\
  5 11 15 17 34 39 40 46 60\\
   8 20 31 36 39 50 51 53 61\\
   4 10 14 16 33 38 39 54 59\\
   5 10 21 24 25 30 37 61 63\\
   4 18 20 23 24 29 45 60 62\\
   3 24 31 35 43 46 48 54 56\\
   8 13 19 24 27 33 40 55 57\\
   1 2 3 4 5 6 7 8 9\\
   1 22 29 33 41 46 52 53 63\\
   1 15 20 21 26 35 42 57 59\\
   1 11 13 16 30 44 45 51 56\\
\hline
\end{tabular}
\caption{Parallel classes of $G_1$}
\label{pclG1}
\end{table}

\begin{table}
\centering
\begin{tabular}{|l|}
\hline
  19 20 21 22 23 24 25 26 27\\
   6 9 10 15 25 30 32 43 44\\
   4 7 13 17 23 28 30 41 42\\
   2 8 12 17 27 32 34 37 45\\
  2 5 11 15 21 28 35 39 40\\
   5 8 14 18 24 29 31 42 43\\
   3 9 13 18 19 33 35 37 38\\
   3 6 12 16 22 29 36 40 41\\
   1 4 10 14 20 34 36 38 39\\
   1 7 11 16 26 31 33 44 45\\
\hline
\end{tabular}
\caption{Parallel classes of $G'_1$}
\label{pclPG1}
\end{table}

\begin{table}
\centering
\begin{tabular}{|l|l|}
\hline
$l_1$ & 14 22 30 35 40 \\
$l_7$ & 12 25 34 38 40 \\
$l_{13}$ &  1 21 28 35 45 \\
$l_{19}$ &  1 13 26 27 43 \\
$l_{25}$ &  1  9 18 30 34 \\
$l_{31}$ &  1  8 23 39 42 \\
$l_{37}$ &  1 11 12 22 32 \\
$l_{43}$ &  1 16 17 36 38 \\
$l_{49}$ &  1 10 15 24 33 \\
$l_{55}$ &  7 10 27 30 44 \\
$l_{58}$ & 13 16 20 23 44 \\
$l_{61}$ & 19 21 23 40 43 \\
$l_{63}$ & 37 38 39 44 45 \\
\hline
\end{tabular}
\caption{Line orbit representatives of $G_2$}
\label{lMr}
\end{table}

\begin{table}
\centering
\begin{tabular}{|l|}
\hline
  1  9 18 30 34 \\
  2 10 13 25 35 \\
  3 11 14 26 36 \\
  4 12 15 27 31 \\
  5  7 16 28 32 \\
  6  8 17 29 33 \\
 19 21 23 40 43 \\
 20 22 24 41 42 \\
 37 38 39 44 45\\
\hline
\end{tabular}
\caption{The unique parallel class of $G_2$}
\label{parlines}
\end{table}

\begin{table}
\centering
\begin{tabular}{|l|}
\hline
   31 32 33 34 35 36 42 43 44 \\
   2 8 15 20 21 26 30 32 38 \\
   5 11 18 23 24 27 29 35 38  \\
   6 12 13 19 24 28 30 36 39  \\
   4 10 17 22 23 26 28 34 37 \\
  8 10 12 14 16 18 41 43 45  \\
   1 7 14 19 20 25 29 31 37  \\
  7 9 11 13 15 17 40 42 45  \\
   1 2 3 4 5 6 40 41 44 \\
  3 9 16 21 22 25 27 33 39 \\
\hline
\end{tabular}
\caption{Parallel classes of $G'_2$}
\label{pclPG2}
\end{table}


\section{Partial geometries arising from planes of small order}
\label{S4}

The upper bounds (\ref{eq6}) and (\ref{eq6'})
of Theorem \ref{t1} apply only to
 sets of parallel classes that are pairwise orthogonal.
It is an interesting open question if the total number of
parallel classes of any partial geometry with parameters (\ref{eq3})
or its dual geometry can exceed any of these bounds.

 Motivated by this question, we examined the parallel classes
of lines in partial geometries arising
from maximal arcs of degree $d$, $2\le d <q$
 in the projective planes of even order $q=2^r \le 16$
via Construction \ref{con}.

The smallest $q$ that meets the condition $2 \le d <q$  
is $q=4$. In this case $d=d'=2$, and any maximal arc is a hyperoval.
The projective plane $PG(2,4)$ of order 4 contains 168 hyperovals, 
all being in one orbit under the collineation group of the plane.
Thus, up to isomorphism, there is a unique partial geometry $pg(2,2,1)$
associated with a hyperoval in $PG(2,4)$, known also as
the generalized quadrangle $W(2)$ (Payne and Thas \cite{PT}).
The total number of parallel classes of lines of $W(2)$
is exactly six, and the six parallel classes are pairwise orthogonal.
We note that the 2-rank of the incidence matrix
of $W(2)$ is 10, and is equal to the 2-rank of 
the incidence matrix of the projective plane of order 4.
 
If $q=8$, the possible degrees of maximal arcs 
are $d=2$ and $d=4$.
All hyperovals in $PG(2,8)$ are projectively equivalent.
Thus, up  to isomorphism, there is a unique partial geometry
$pg(6,4,3)$ arising from a hyperoval in $PG(2,8)$, 
being isomorphic to the dual of the partial geometry 
$G_1$ from Section \ref{S3}, while
every partial geometry $pg(4,6,3)$ arising from a
maximal arc of degree 4 in $PG(2,8)$ is isomorphic to 
 $G_1$. It was shown in  Section \ref{S3} that
$G_1$ has a total of exactly 28 parallel classes of lines,
every two being orthogonal, and its dual geometry has exactly 10
parallel classes of lines, every two being orthogonal.
As we mentioned in Section \ref{S3},
the 2-rank of the incidence matrix
of $G_1$ is 28, which is also the 2-rank 
of the incidence matrix
of the projective plane of order 8.

If $q=16$, the possible  degrees of a maximal arc are
 $d=2$, $d=4$, and $d=8$.
If $d=2$, a maximal arc of degree 2,
that is, a hyperoval $H$, gives rise to a partial geometry
$pg(14,8,7)$, while the dual arc of $H$ is of degree 8 and
gives rise to
the dual partial geometry with parameters $pg(8,14,7)$.
If $d=4$, a maximal arc of degree 4 and its
dual arc of degree 4 give rise to two not
necessarily isomorphic partial geometries
with parameters $pg(12,12,9)$.

It follows from Theorem \ref{t2}
that a partial geometry $pg(14,8,7)$ associated with
a hyperoval in a projective plane of order 16 must have a set
of 18 pairwise orthogonal parallel classes of lines, and its
dual geometry must have 120 pairwise orthogonal
parallel classes of lines, while any $pg(12,12,9)$ arising
from a maximal arc of degree 4, as well as its dual geometry each
must have a set of 52 pairwise orthogonal parallel classes.

Up to duality, there are 22 nonisomorphic projective planes
of order 16 that are known currently: four planes are self-dual
and nine planes are not self-dual (see \cite{Moor}, \cite{PRS}).
All hyperovals in the known 22 projective planes of order
16 were enumerated and classified up to equivalence\footnote{Two
hyperovals, or more generally, 
two maximal arcs in a projective plane $\cal{P}$ are equivalent if  
one can be obtained from the other by applying an automorphism of 
$\cal{P}$.}
by  Penttila,  Royle, and Simpson \cite{PRS}:
altogether, there are 93 equivalence classes of hyperovals
in the 22 planes, and by duality, 93 equivalence classes of
maximal arcs of degree 8.
The specific line sets of the known projective planes of order 16
and representatives of the equivalence classes of hyperovals
 were graciously provided 
to the second author by Gordon F. Royle.

Using Magma \cite{magma} and Cliquer \cite{cliquer},
we checked by computer that the 93 equivalence classes of
hyperovals give rise to 
93 nonisomorphic 
partial geometries $pg(14,8,7)$ arising from hyperovals,
and, by duality, 93 nonisomorphic $pg(8,14,7)$ arising
from maximal arcs of degree 8 in projective planes of order 16.
We computed  all parallel classes of lines,
and in all cases, a partial geometry
$pg(14,8,7)$ associated with a hyperoval has
exactly 18 parallel classes, while its dual geometry
has exactly 120 parallel classes: note that 18 and 120 are 
the upper bounds  (\ref{eq6}) and (\ref{eq6'})
from Theorem \ref{t1} respectively. 
In addition, the 2-rank of the incidence matrix
of each partial geometry $pg(14,8,7)$  is equal to the 2-rank of the 
underlying projective plane.

Table \ref{tab7} contains data about the partial geometries $pg(8,14,7)$ 
arising from maximal arcs of degree 8.

\begin{Note}
\label{n2}
{\rm
An interesting {\bf open problem} is to find an analogue of 
Mathon's partial geometry
when $q=16$, that is, a partial geometry $pg(14,8,7)$
that does not arises from any hyperoval  
in a projective plane of order 16, 
or equivalently, according to Theorem \ref{t2}, 
its dual geometry does not contain
a set of 120 pairwise orthogonal parallel classes of lines.
The line graph of such a partial geometry has to be isomorphic to
$\bar{T}(18)$, and every line  corresponds to a clique of size
nine of $\bar{T}(18)$. Since the graph  $\bar{T}(18)$ contains
\[
\frac{18!}{2^{9}9!} = 34,459,425
\]
cliques of size 9, finding a collection of 255 cliques of size 9,
every two sharing at most one vertex, 
would be a very challenging computational problem, at the least.
}
\end{Note}

If $d=4$, the only projective plane of order 16 for which
all maximal arcs of degree 4 have been classified
up to equivalence, is the Desarguesian plane $PG(2,16)$
(Ball and Blokhuis \cite{BB}).
The maximal arcs of degree 4 have not been classified completely
in any of the known  non-Desarguesian planes
of order 16,  although such arcs have been found in
all but four of the known non-Desarguesian planes of order 16
\cite{GMT}, \cite{GTW},  \cite{HST}.
All currently known maximal arcs of degree 4
 are available online at
\begin{verbatim} http://pages.mtu.edu/~tonchev/pointsetsOFmaxArcs.txt,
\end{verbatim}
where the maximal arcs are listed as sets of points in the
the projective planes of order 16 
that are available online at
 \begin{verbatim}
 http://pages.mtu.edu/~tonchev/planesOForder16.txt
 \end{verbatim}

Using Magma \cite{magma}, we classified, up to isomorphism,
the partial geometries $pg(12,12,9)$ arising from the known
maximal arcs of degree 4 in projective planes of order 16,
and computed the 2-rank of their incidence matrices
and all parallel classes of lines.
Table \ref{tab8} summarizes the data about these partial geometries.
Up to isomorphism, there are 59 such geometries (13 self-dual
and 23 non-self-dual ones). Similarly to the previous cases 
for degrees other than 4, or planes of smaller order,
the 2-rank of each partial geometry turned out to be equal to the 2-rank
of the corresponding plane, and the total number of parallel 
classes always  
matches the upper bound  of Theorem \ref{t1}.

These observations motivate us to formulate the following conjectures. 
 
\begin{Con}
\label{conj1}
Let $G$ be a partial geometry with parameters (\ref{eq3})
for some integer numbers $d\ge 2$, $d'\ge 2$.
Then either $G$ does not arise from 
a maximal arc of degree $d$ in a projective plane
of order $q=dd'$, or $G$ is obtainable from exactly
one projective plane of order $q=dd'$ via Construction \ref{con}.
\end{Con}
\begin{Con}
\label{conj2}
If a partial geometry $G$ with parameters (\ref{eq3})
for some $d=2^i \ge 2$, $d'= 2^j \ge 2$
arises from a maximal arc of degree $d$
in a projective plane $\cal{P}$ of order $q=2^{i+j}$
via Construction \ref{con}, then the 2-rank
of the $(0,1)$-incidence matrix of $G$ is equal to the 2-rank of 
the $(0,1)$-incidence matrix of $\cal{P}$.
\end{Con}



\newpage

\section{Appendix}

\begin{table}[htb!]
\centering
\scalebox{.7}{
\begin{tabular}{|c|c|c|c|c|c|}
  \hline
No. & Hyperoval & $|Aut(G)|$ & 2-rank & \# Par. cl.  \\
  \hline
1 & PG(2,16)hyp.1 & 144 & 82 & 120 \\
  \hline
2 & PG(2,16)hyp.2 & 16320 & 82 & 120  \\
  \hline
3 & DEMPhyp.1 & 16 & 102 & 120  \\
\hline
4 & DEMPhyp.2 & 16 & 101 & 120   \\
  \hline
5 & DEMPhyp.3 & 16 & 102 & 120   \\
  \hline
6 & DEMPhyp.4 & 16 & 102 & 120   \\
  \hline
7 & DEMPhyp.5 & 16 & 100 & 120   \\
  \hline
8 & DEMPhyp.6 & 16 & 102 & 120   \\
  \hline
9 & DEMPhyp.7 & 16 & 101 & 120   \\
  \hline
10 & DEMPhyp.8 & 16 & 102 & 120  \\
  \hline
11 & DEMPhyp.9 & 16 & 101 & 120   \\
  \hline
12 & DEMPhyp.10 & 16 & 102 & 120  \\
  \hline
13 & DEMPhyp.11 & 16 & 101 & 120  \\
  \hline
14 & DEMPhyp.12 & 64 & 100 & 120  \\
  \hline
15 & DEMPhyp.13 & 64 & 100 & 120  \\
  \hline
16 & DEMPhyp.14 & 64 & 100 & 120  \\
  \hline
17 & DEMPhyp.15 & 80 & 102 & 120  \\
  \hline
18 & dDEMPhyp.1 & 2 & 102 & 120  \\
  \hline
19 & dDEMPhyp.2 & 6 & 102 & 120  \\
  \hline
20 & SEMI4hyp.1 & 16 & 98 & 120  \\
  \hline
21 & SEMI4hyp.2 & 16 & 98 & 120  \\
  \hline
22 & SEMI4hyp.3 & 16 & 98 & 120  \\
  \hline
23 & SEMI2hyp.1 & 3 & 98 & 120   \\
  \hline
24 & SEMI2hyp.2 & 8 & 98 & 120   \\
  \hline
25 & SEMI2hyp.3 & 8 & 98 & 120   \\
  \hline
26 & SEMI2hyp.4 & 16 & 98 & 120  \\
  \hline
27 & SEMI2hyp.5 & 16 & 98 & 120  \\
  \hline
28 & SEMI2hyp.6 & 16 & 97 & 120  \\
  \hline
29 & SEMI2hyp.7 & 16 & 98 & 120  \\
  \hline
30 & SEMI2hyp.8 & 16 & 97 & 120  \\
  \hline
31 & SEMI2hyp.9 & 16 & 98 & 120  \\
  \hline
32 & SEMI2hyp.10 & 16 & 97 & 120 \\
  \hline
33 & SEMI2hyp.11 & 16 & 98 & 120 \\
  \hline
34 & SEMI2hyp.12 & 16 & 98 & 120 \\
  \hline
35 & SEMI2hyp.13 & 16 & 97 & 120 \\
  \hline
36 & SEMI2hyp.14 & 26 & 98 & 120 \\
  \hline
37 & SEMI2hyp.15 & 16 & 98 & 120 \\
  \hline
38 & SEMI2hyp.16 & 32 & 98 & 120 \\
  \hline
39 & SEMI2hyp.17 & 32 & 97 & 120 \\
  \hline
40 & LMRHhyp.1 & 16 & 100 & 120 \\
 \hline
41 & LMRHhyp.2 & 16 & 103 & 120 \\
 \hline
42 & LMRHhyp.3 & 16 & 103 & 120 \\
 \hline
43 & LMRHhyp.4 & 16 & 103 & 120 \\
 \hline
44 & LMRHhyp.5 & 64 & 103 & 120 \\
 \hline
45 & LMRHhyp.6 & 112 & 99 & 120 \\
 \hline
46 & dLMRHhyp.1 & 14 & 105 & 120 \\
\hline
\end{tabular}}
\caption{The pg(8,14,7)'s arising from maximal (120,8)-arcs}
\label{tab7}
\end{table}

\begin{table}[htb!]
\ContinuedFloat  
\centering
\scalebox{.70}{
\begin{tabular}{|c|c|c|c|c|c|}
  \hline
No. & Hyperoval & $|Aut(G)|$ & 2-rank & \# Par. cl.  \\
  \hline
47 & MATHhyp.1 & 8 & 107 & 120   \\
 \hline
48 & dMATHhyp.1 & 4 & 108 & 120  \\
 \hline
49 & dMATHhyp.2 & 4 & 108 & 120  \\
 \hline
50 & dMATHhyp.3 & 8 & 109 & 120 \\
 \hline
51 & HALLhyp.1 & 16 & 97 & 120   \\
  \hline
52 & HALLhyp.2 & 64 & 98 & 120  \\
  \hline
53 & HALLhyp.3 & 64 & 97 & 120  \\
  \hline
54 & HALLhyp.4 & 320 & 97 & 120 \\
  \hline
55 & dHALLhyp.1 & 2 & 98 & 120 \\
  \hline
56 & dHALLhyp.2 & 2 & 98 & 120  \\
  \hline
57 & dHALLhyp.3 & 6 & 98 & 120  \\
  \hline
58 & BBH1hyp.1 & 8 & 107 & 120  \\
  \hline
59 & BBH1hyp.2 & 16 & 109 & 120  \\
\hline
60 & BBH1hyp.3 & 32 & 107 & 120 \\
\hline
61 & JOWKhyp.1 & 16 & 99 & 120 \\
  \hline
62 & JOWKhyp.2 & 16 & 99 & 120  \\
  \hline
63 & JOWKhyp.3 & 16 & 100 & 120 \\
  \hline
64 & JOWKhyp.4 & 16 & 99 & 120  \\
  \hline
65 & JOWKhyp.5 & 64 & 99 & 120  \\
  \hline
66 & JOWKhyp.6 & 112 & 99 & 120 \\
  \hline
67 & dJOWKhyp.1 & 14 & 99 & 120 \\
  \hline
68 & JOHNhyp.1 & 16 & 111 & 120 \\
  \hline
69 & DSFPhyp.1 & 16 & 103 & 120 \\
  \hline
70 & DSFPhyp.2 & 16 & 103 & 120 \\
\hline
71 & DSFPhyp.3 & 16 & 103 & 120 \\
  \hline
72 & DSFPhyp.4 & 16 & 103 & 120 \\
  \hline
73 & DSFPhyp.5 & 16 & 103 & 120 \\
  \hline
74 & DSFPhyp.6 & 16 & 103 & 120 \\
  \hline
75 & DSFPhyp.7 & 16 & 103 & 120 \\
  \hline
76 & DSFPhyp.8 & 16 & 103 & 120 \\
  \hline
77 & DSFPhyp.9 & 16 & 103 & 120 \\
  \hline
78 & DSFPhyp.10 & 16 & 103 & 120 \\
  \hline
79 & DSFPhyp.11 & 16 & 103 & 120 \\
  \hline
80 & DSFPhyp.12 & 16 & 103 & 120 \\
  \hline
81 & DSFPhyp.13 & 16 & 103 & 120 \\
  \hline
82 & DSFPhyp.14 & 16 & 103 & 120 \\
  \hline
83 & DSFPhyp.15 & 16 & 103 & 120 \\
  \hline
84 & DSFPhyp.16 & 16 & 102 & 120 \\
  \hline
85 & DSFPhyp.17 & 16 & 103 & 120 \\
  \hline
86 & DSFPhyp.18 & 16 & 103 & 120 \\
  \hline
87 & DSFPhyp.19 & 16 & 103 & 120 \\
  \hline
88 & DSFPhyp.20 & 16 & 103 & 120 \\
  \hline
89 & DSFPhyp.21 & 16 & 103 & 120 \\
  \hline
90 & DSFPhyp.22 & 64 & 103 & 120 \\
  \hline
91 & BBH2hyp.1 & 4 & 111 & 120  \\
  \hline
92 & BBH2hyp.2 & 4 & 111 & 120  \\
  \hline
93 & BBS4hyp.1 & 16 & 111 & 120 \\
  \hline
\end{tabular}}
\caption{The pg(8,14,7)'s arising from maximal (120,8)-arcs (continued).}
\label{tab7}
\end{table}

\newpage

 \begin{table}
\centering
\scalebox{.8}{
\begin{tabular}{|c|c|c|c|c|c|c|}
  \hline
 & Maximal & & 2-rank & $\#$ & Isomorphic & Isomorphic  \\
$\#$ & 52-arc & $|Aut(G)|$ &  & of & to & to \\
 &  &   &  & par. clas. & its dual? & others?  \\

\hline
1 & PG(2,16).1 & 68 & 82 & 52 & Yes & No \\
  \hline
2 & PG(2,16).2 & 408 &  82 & 52 & Yes & No \\
  \hline
3 & DEMP.1 & 24 &  102 & 52 & No & No \\
  \hline
4 & DEMP.2 & 144 & 102 & 52 & No & No \\
  \hline
5 & DEMP.3 & 24 & 102 & 52 & No & No \\
  \hline
6 & DEMP.4 & 48 & 102 & 52 & No & No \\
  \hline
7 & DEMP.5 & 4 & 102 & 52 & No & No \\
  \hline
8 & SEMI4.1 & 96 & 98 & 52 & Yes & No \\
  \hline
9 & SEMI2.1 & 24 & 98 & 52 & Yes & No \\
  \hline
10 & SEMI2.2 & 144 & 98 & 52 & Yes & No \\
  \hline
11 & SEMI2.3 & 32 & 98 & 52 & Yes & No \\
  \hline
12 & SEMI2.4 & 32 & 98 & 52 & Yes & No \\
  \hline
13 & SEMI2.5 & 16 & 98 & 52 & Yes & No \\
  \hline
14 & SEMI2.6 & 48 & 98 & 52 & No & SEMI2.7 \\
  \hline
15 & SEMI2.7 & 48 & 98 & 52 & No & SEMI2.6 \\
  \hline
16 & LMRH.1 & 96 & 106 & 52 & No & No \\
  \hline
17 & LMRH.2 & 32 & 106 & 52 & No & No \\
  \hline
18 & MATH.1 & 24 & 109 & 52 & No & No \\
  \hline
19 & MATH.2 & 32 & 108 & 52 & No & No \\
  \hline
20 & MATH.3 & 32 & 108 & 52 & No & No \\
  \hline
21 & MATH.4 & 32 & 108 & 52 & No & No \\
  \hline
22 & MATH.5 & 16 & 109 & 52 & No & No \\
  \hline
23 & MATH.6 & 16 & 109 & 52 & No & No \\
  \hline
24 & MATH.7 & 16 & 109 & 52 & No & No \\
  \hline
25 & HALL.1 & 24 & 98 & 52 & No & No \\
  \hline
26 & HALL.2 & 4 & 98 & 52 & No & No \\
  \hline
27 & BBH1.1 & 24 & 110 & 52 & Yes & No \\
  \hline
28 & BBH1.2 & 32 & 110 & 52 & Yes & No \\
  \hline
29 & BBH1.3 & 4 & 110 & 52 & Yes & No \\
  \hline
30 & JOWK.1 & 16 & 100 & 52 & No & No \\
  \hline
31 & JOWK.2 & 32 & 100 & 52 & No & No \\
  \hline
32 & JOHN.1 & 32 & 113 & 52 & No & No \\
  \hline
33 & JOHN.2 & 32 & 113 & 52 & No & No \\
  \hline
34 & JOHN.3 & 32 & 113 & 52 & No & No \\
  \hline
35 & JOHN.4 & 32 & 113 & 52 & No & No \\
  \hline
36 & DSFP.1 & 24 & 106 & 52 & No & No \\
  \hline
\end{tabular}}
\caption{The pg(12,12,9)'s arising from maximal (52,4)-arcs}
\label{tab8}
\end{table}

\end{document}